\documentclass[11pt]{amsart} 
\usepackage[T1]{fontenc}
\usepackage[latin1]{inputenc}
\usepackage{amsmath}
\usepackage{amssymb}
\usepackage{graphics}
\usepackage{epic}
\newtheorem{theorem}{Theorem}

\begin{document}
\title{The cover pebbling theorem}
\author[J.~Sjöstrand]{Jonas~Sjöstrand}
\address{Dept. of Mathematics, KTH, SE-100 44 Stockholm, Sweden}
\email{jonass@kth.se}
\date{October 6, 2004}

\begin{abstract}
For any configuration of pebbles on the nodes of a graph,
a pebbling move replaces two pebbles on one node by one
pebble on an adjacent node. A cover pebbling is a move 
sequence ending with no empty nodes. The number of pebbles 
needed for a cover pebbling starting with all pebbles on 
one node is trivial to compute and it was conjectured
that the maximum of these simple cover pebbling numbers
is indeed the general cover pebbling number of the graph.
That is, for any configuration of this size, there exists
a cover pebbling. In this note, we prove a generalization of the conjecture.
All previously published results about cover pebbling numbers
for special graphs (trees, hypercubes etcetera) are
direct consequences of this theorem. We also prove that
the cover pebbling number of a product of two graphs equals
the product of the cover pebbling numbers of the graphs.
\end{abstract}
\maketitle

\subsection*{Introduction}\noindent
Pebbling, peg solitaire, chip firing and checker jumping
are some kindred combinatorial games on graphs. Put some 
tokens on the nodes, define local moves and you can start
asking questions about convergence, reachability and 
enumeration!
But this is not a true description of how these games
came into existence. Each one has its own roots in areas
such as number theory, statistical mechanics, economics
and of course recreational mathematics.

The pebbling game appeared in the 1980s and
comes in two flavours. In the first version, played on
a directed graph, a move consists in replacing a pebble
on one node by new pebbles on the adjacent nodes, 
moving along directed edges. The 1995 paper by Eriksson
\cite{eriksson} seems to have solved this game
completely. 

The second pebbling version, which is still very hot,
was introduced in 1989 by Chung \cite{chung} and is
played on a connected graph, directed or undirected.
A move replaces two pebbles on one node by one pebble
on an adjacent node, and this is the pebbling rule for
the remainder of this paper.

The most important reachability questions concern
the {\em pebbling number} and the {\em cover pebbling 
number} of a graph, that is the smallest
$n$ such that from any initial distribution of $n$ 
pebbles, it is possible to pebble any
desired node respectively pebble all nodes. In a series of
recent papers by Crull et al.~\cite{crull}, Watson and Yerger \cite{watson},
Hurlbert and Munyan \cite{hurlbert2}, and
Tomova and Wyels \cite{tomova}, 
the cover pebbling number has been derived for
several classes of graphs. These results are all special
cases of our main theorem, conjectured by Crull et al.~in
\cite{crull}.

\subsection*{General covers and simple distributions.}\noindent
Following Crull et al.~we generalize the situation like this:
Instead of trying to place at least one pebble on each node,
we define a goal distribution $w$ of pebbles.
A {\em $w$-cover} is a distribution of pebbles such that
every node has at least as many pebbles as in $w$. We
write $w(v)$ for the number of pebbles on the node $v$ in $w$.
In this terminology,
the usual cover is the special case where $w$ is the {\em 1-distribution},
i.e.~$w(v)=1$ for all nodes $v$. The {\em $w$-cover pebbling number}
is the smallest $n$ such that, from any initial
distribution of $n$ pebbles, it is possible to obtain a $w$-cover.
We will require $w$ to be positive, i.e.~there should
be at least one pebble on each node.

A node $v$ is {\em fat}, {\em thin} respectively
{\em perfect} if the number of pebbles on it is greater than, less than,
respectively equal to $w(v)$.

The initial distribution is said to be {\em simple}
if all pebbles are on one single node. For two nodes $v$ and $u$,
the {\em distance
$d(v,u)$ from $v$ to $u$} is the length of the minimal
path from $v$ to $u$. (For a directed graph, $d(v,u)\neq d(u,v)$
in general.) The {\em cost} from a node $v$
of a pebble on a node $u$ is $2^{d(v,u)}$,
and the sum of the costs from $v$ of all pebbles in $w$
is the cost of cover pebbling from $v$.

In the graph
below, $8+8+4+2+1$ pebbles on $v$ are necessary and
sufficient for a cover pebbling if $w$ is the 1-distribution.

\begin{figure}[h]
\begin{picture}(70,36)(-18,-18)
\put(40,0){\vector(-1,0){17}}
\put(20,0){\vector(-1,0){17}}
\put(0,0){\vector(-1,1){12}}
\put(0,0){\vector(-1,-1){12}}
\put(-11,14){\vector(4,-1){48}}
\put(-11,-14){\vector(4,1){48}}
\put(-14,14){\circle*{7}}
\put(-14,-14){\circle*{7}}
\put(0,0){\circle*{7}}
\put(20,0){\circle*{7}}
\put(40,0){\circle*{7}}
\put(45,-1){$v$}
\end{picture}
\caption{A cover pebbling from $v$ needs 23 pebbles.}
\label{fig:graph}
\end{figure}
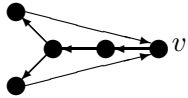\noindent
In each pebbling move, the total number of pebbles decreases, but
the total {\em value} is invariant, if the value of a pebble is defined
to be the number of pebbles that have gone into it. Recursively speaking,
the value of a newborn pebble is the sum of the values of its demised
parents, the original pebbles being unit valued.

\subsection*{The cover pebbling theorem}\noindent
For nonsimple initial distributions, costs are ill-defined and
there is no easy way to see if a cover pebbling exists.
The following theorem tells us not to worry about that when it
comes to computing the cover pebbling number of a graph, for
this number is always determined by a simple distribution.    

\begin{theorem}
Let $w$ be a positive goal distribution.
To determine the $w$-cover pebbling number of a (directed or undirected)
connected graph, it is sufficient to
consider simple initial distributions.
In fact, for any initial distribution that admits no cover pebbling,
all pebbles may be concentrated to one of the fat nodes\footnote{
Of course, this is not true if there are no fat nodes, but then
any node will do.
} with cover pebbling still not possible.
\end{theorem}
\begin{proof}
Start with a distribution that admits no cover pebbling.

If there are no fat nodes, we can concentrate all the
pebbles to any of the nodes. The cost of cover pebbling from
this node is of course no less than the number of pebbles in $w$,
so cover pebbling is still not possible.

If some node is fat, we will have to do some pebbling.
During the pebbling we will always maintain the following
efficiency condition:
{\em Every pebble has a value no greater
than the cost from its nearest fat node (the fat node that
minimizes this cost).}
At the beginning all pebbles have the value one, so the efficiency
condition is trivially satisfied.

Now pebble like this:
Among all pairs $(f,t)$ of a fat and a thin node, take
one that minimizes the distance $d(f,t)$.
Let $fp_1p_2\cdots p_{d-1}t$ be a minimal path from $f$ to $t$.
Every inner node $p_i$
of this path must be perfect, since if it were thin,
then $(f,p_i)$ would be a (fat,thin)-pair with $d(f,p_i)<d(f,t)$,
and if it were fat, then $(p_i,t)$ would be a (fat,thin)-pair
with $d(p_i,t)<d(f,t)$. Furthermore, $f$ must
be a nearest fat node to $t$ and to every $p_i$. Now we play two pebbles on
$f$ to $p_1$, then we play the new pebble on $p_1$ together with any old
pebble from $p_1$ to $p_2$, then the new pebble and an old one on $p_2$ to
$p_3$, and so on, until we reach $t$.

\begin{figure}[h]
\begin{picture}(130,30)(-2,-15)
\put(0,0){\vector(1,0){17}}
\put(20,0){\vector(1,0){17}}
\put(40,0){\vector(1,0){17}}
\dashline[-7]{2}(60,0)(100,0)
\put(80,0){\vector(1,0){17}}
\put(100,0){\vector(1,0){17}}
\put(0,0){\circle*{7}}
\put(20,0){\circle*{7}}
\put(40,0){\circle*{7}}
\put(100,0){\circle*{7}}
\put(120,0){\circle*{7}}
\put(-2,10){$f$}
\put(16,10){$p_1$}
\put(36,10){$p_2$}
\put(90,10){$p_{d-1}$}
\put(118,10){$t$}
\put(-5,-15){$17$}
\put(18,-15){$1$}
\put(38,-15){$1$}
\put(98,-15){$1$}
\put(118,-15){$0$}
\end{picture}
\\
\begin{picture}(130,40)(-2,-25)
\put(60,-10){\Huge $\Downarrow$}
\end{picture}
\\
\begin{picture}(130,30)(-2,-15)
\put(0,0){\vector(1,0){17}}
\put(20,0){\vector(1,0){17}}
\put(40,0){\vector(1,0){17}}
\dashline[-7]{2}(60,0)(100,0)
\put(80,0){\vector(1,0){17}}
\put(100,0){\vector(1,0){17}}
\put(0,0){\circle*{7}}
\put(20,0){\circle*{7}}
\put(40,0){\circle*{7}}
\put(100,0){\circle*{7}}
\put(120,0){\circle*{7}}
\put(-2,10){$f$}
\put(16,10){$p_1$}
\put(36,10){$p_2$}
\put(90,10){$p_{d-1}$}
\put(118,10){$t$}
\put(-5,-15){$15$}
\put(18,-15){$0$}
\put(38,-15){$0$}
\put(98,-15){$0$}
\put(118,-15){$1$}
\end{picture}
\caption{Playing two pebbles from $f$ and continuing all
the way to $t$ in the case where $w$ is the 1-distribution.}
\end{figure}
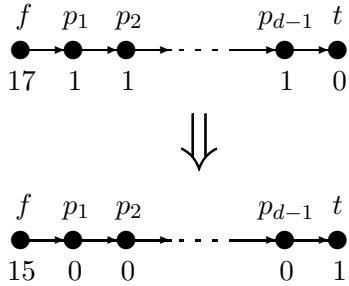\noindent
The value of the new pebble on $t$ is 2 plus the sum of the values
of the old pebbles on $p_1,\ldots,p_{d-1}$ that were consumed.
By the efficiency condition
this is no greater than $2+2^1+2^2+\cdots+2^{d-1}$ which equals $2^d$.
Thus the condition is satisfied even after this operation.
It is possible that $f$ is no longer fat, but this only makes
it easier to fulfil the condition.

We iterate the above procedure (choosing a new pair $(f,t)$, and so on)
until no node is fat.
During each iteration the total number of pebbles on fat nodes decreases,
so we cannot continue forever.

Let $f$ be the fat node that survived the longest.
Then each pebble value is at most equal to its cost from $f$.
But there are still thin nodes, so the cost of cover pebbling from $f$
exceeds the total value of the pebbles. Therefore,
cover pebbling is not possible with all pebbles initially on $f$.
\end{proof}

\subsection*{The cover number for some classes of graphs}
\noindent
The cover number is now easy to compute for any graph.
Here is a table of some classes of undirected graphs, for
the case that $w$ is the 1-distribution.

\noindent
\begin{center}
\begin{tabular}{|l|l@{\hspace{-3mm}}c@{\hspace{0mm}}|c|}
\hline {\bf class} & \multicolumn{2}{c|}{\bf example} & {\bf cover number} \\

\hline\hline $n$-path &
$n=4$: &
\begin{tabular}{c}
\begin{picture}(70,10)(-5,-5)
\put(0,0){\line(1,0){60}}
\put(0,0){\circle*{7}}
\put(20,0){\circle*{7}}
\put(40,0){\circle*{7}}
\put(60,0){\circle*{7}}
\end{picture}
\end{tabular}
& $2^n-1$ \\

\hline $2n$-cycle &
$n=2$: &
\begin{tabular}{c}
\begin{picture}(30,30)(-5,-5)
\put(0,0){\line(1,0){20}}
\put(0,0){\line(0,1){20}}
\put(0,20){\line(1,0){20}}
\put(20,0){\line(0,1){20}}
\put(0,0){\circle*{7}}
\put(20,0){\circle*{7}}
\put(0,20){\circle*{7}}
\put(20,20){\circle*{7}}
\end{picture}\end{tabular}
& $3\cdot(2^n-1)$ \\

\hline $(2n-1)$-cycle &
$n=2$: &
\begin{tabular}{c}
\begin{picture}(30,30)(-5,-5)
\put(0,0){\line(1,0){20}}
\put(0,0){\line(3,5){10}}
\put(10,16.7){\line(3,-5){10}}
\put(0,0){\circle*{7}}
\put(20,0){\circle*{7}}
\put(10,16.7){\circle*{7}}
\end{picture}
\end{tabular}
& $2^{n+1}-3$ \\

\hline\begin{tabular}{@{\hspace{0mm}}l}
$n$-dimensional \\
hypercube\end{tabular} &
$n=3$: &
\begin{tabular}{c}
\begin{picture}(40,40)(-5,-5)
\put(0,0){\circle*{7}}
\put(20,0){\circle*{7}}
\put(0,20){\circle*{7}}
\put(20,20){\circle*{7}}
\put(10,10){\circle*{7}}
\put(30,10){\circle*{7}}
\put(10,30){\circle*{7}}
\put(30,30){\circle*{7}}

\put(0,0){\line(1,0){20}}
\put(0,0){\line(0,1){20}}
\put(0,20){\line(1,0){20}}
\put(20,0){\line(0,1){20}}

\put(10,10){\line(1,0){20}}
\put(10,10){\line(0,1){20}}
\put(10,30){\line(1,0){20}}
\put(30,10){\line(0,1){20}}

\put(0,0){\line(1,1){10}}
\put(20,0){\line(1,1){10}}
\put(0,20){\line(1,1){10}}
\put(20,20){\line(1,1){10}}
\end{picture}
\end{tabular}
& $3^n$ \\

\hline complete graph $K_n$ &
$n=4$: &
\begin{tabular}{c}
\begin{picture}(30,30)(-5,-5)
\put(0,0){\line(1,0){20}}
\put(0,0){\line(3,5){10}}
\put(10,16.7){\line(3,-5){10}}
\put(0,0){\line(4,1){26.7}}
\put(18,-2){\line(1,1){10}}
\put(10,16.7){\line(5,-3){16.7}}

\put(0,0){\circle*{7}}
\put(20,0){\circle*{7}}
\put(10,16.7){\circle*{7}}
\put(26.7,6.7){\circle*{7}}
\end{picture}
\end{tabular}
& $2n-1$ \\

\hline \begin{tabular}{@{\hspace{0mm}}l}
complete multipartite\\
graph $K_{n_1,\ldots,n_k}$ where\\
$n_1\geq\cdots\geq n_k$\end{tabular} &
\begin{tabular}{@{\hspace{0mm}}l}
$n_1=3$\\
$n_2=2$
\end{tabular} &
\begin{tabular}{c}
\begin{picture}(50,50)(-5,-25)
\put(0,0){\circle*{7}}
\put(0,20){\circle*{7}}
\put(0,-20){\circle*{7}}
\put(40,10){\circle*{7}}
\put(40,-10){\circle*{7}}

\put(0,0){\line(4,1){40}}
\put(0,0){\line(4,-1){40}}
\put(0,20){\line(4,-1){40}}
\put(0,20){\line(4,-3){40}}
\put(0,-20){\line(4,1){40}}
\put(0,-20){\line(4,3){40}}
\end{picture}\end{tabular}
&
$4n_1+2n_2+\cdots+2n_k-3$
\\

\hline $n$-wheel &
$n=4$: &
\begin{tabular}{c}
\begin{picture}(30,30)(-5,-5)
\put(0,0){\line(1,0){20}}
\put(0,0){\line(0,1){20}}
\put(0,20){\line(1,0){20}}
\put(20,0){\line(0,1){20}}

\put(0,0){\line(1,1){20}}
\put(0,20){\line(1,-1){20}}

\put(0,0){\circle*{7}}
\put(20,0){\circle*{7}}
\put(0,20){\circle*{7}}
\put(20,20){\circle*{7}}
\put(10,10){\circle*{7}}
\end{picture}
\end{tabular}
& $4n-5$ \\ \hline
\end{tabular}
\end{center}

\noindent
All the results in the table were previously known,
but now, in the light of our theorem, they are reduced to simple
exercises.

\subsection*{Product graphs}
One of the conjectures in \cite{hurlbert2} is also an easy consequence
of our theorem. We will prove it in a much more general form.
Introduce the notation $G_1\Box G_2$ for the
product\footnote{
$V(G_1\Box G_2)=V(G_1)\times V(G_2)$ and there is an
edge from $(u_1,u_2)$ to $(v_1,v_2)$ if
$u_1=v_1$ and there is an edge from $u_2$ to $v_2$ in $G_2$,
or if $u_2=v_2$ and there is an edge from $u_1$ to $v_1$ in $G_1$.}
of two (directed or undirected) graphs.

Let $w_1$ and $w_2$ be goal distributions on $G_1$ respectively
$G_2$. Define a goal distribution $w_1\Box w_2$ on $G_1\Box G_2$
by $(w_1\Box w_2)(v_1,v_2)=w_1(v_1)w_2(v_2)$.
Finally, let $\gamma_w$ denote the $w$-cover pebbling number.
Then we have the following theorem.
\begin{theorem}
$$\gamma_{w_1\Box w_2}(G_1\Box G_2)=\gamma_{w_1}(G_1)\gamma_{w_2}(G_2).$$
\end{theorem}
\begin{proof}
The distance from $(v_1,v_2)$ to $(u_1,u_2)$ in $G_1\Box G_2$ is
equal to the sum of the distance from $v_1$ to $u_1$ in $G_1$
and the distance from $v_2$ to $u_2$ in $G_2$, so, for any $v_1$ and $v_2$,
$$\gamma_{w_1\Box w_2}(G_1\Box G_2)
\geq\sum_{(u_1,u_2)\in V(G_1\Box G_2)}
(w_1\Box w_2)(u_1,u_2)2^{d((v_1,v_2),(u_1,u_2))}=$$
$$=\sum_{u_1\in V(G_1)}w_1(u_1)2^{d(v_1,u_1)}
\sum_{u_2\in V(G_2)}w_2(u_2)2^{d(v_2,u_2)}
\leq\gamma_{w_1}(G_1)\gamma_{w_2}(G_2).$$
By the cover pebbling theorem, we can choose $(v_1,v_2)$ to
make the first inequality an equality.
On the other hand, there are $v_1$ and $v_2$ that make the second
inequality an equality.
\end{proof}

\subsection*{Acknowledgement}\noindent
I would like to thank Henrik Eriksson for introducing me to the problem,
and Elin Svensson for giving helpful comments.

\end{document}